\documentclass[11pt,a4paper,reqno]{amsart}  
\usepackage{amssymb} \usepackage{amscd} \usepackage{times}
\newcommand{\bea} {\begin{eqnarray}}
\newcommand{\eea} {\end{eqnarray}}
\newcommand{\Bea} {\begin{eqnarray*}}
\newcommand{\Eea} {\end{eqnarray*}}

\def\zbb{\mathbb{Z}}  
  
  \def\phi{\varphi}
 \def\p1{{\mathbb{P}^1_\zbb}}

\newtheorem{Theorem}{\quad Theorem}[section]

\usepackage{amsmath,  amsfonts}

\usepackage[T1]{fontenc}

\newcommand{\be} {\begin{equation}}
\newcommand{\ee} {\end{equation}}

\begin{document}
\title {Uniform result for solutions of an equation with boundary singularity.}
\author{Samy Skander Bahoura} 
\address{Departement de Mathematiques, Universite Pierre et Marie Curie, 2 place Jussieu, 75005, Paris, France.}
\email{samybahoura@yahoo.fr, samybahoura@gmail.com}
\maketitle
\begin{abstract}

We consider a variational problem with boundary singularity and Dirichlet condition. We give blow-up analysis for sequences of solutions of an equation with exponential nonlinearity. Also we derive a compactness criterion for this problem under some condition. 

\end{abstract}

{ \small  Keywords: blow-up, boundary, singularity, a priori estimate, analytic domain, Lipschitz condition.}

\section{Introduction and Main Results} 

We set $ \Delta = \partial_{11} + \partial_{22} $  on open analytic domain $ \Omega $ of $ {\mathbb R}^2 $.

\bigskip

We consider the following equation:

$$ (P)   \left \{ \begin {split} 
      -\Delta u & = |x|^{2\beta}V e^{u} \,\, &\text{in} \,\, & \Omega  \subset {\mathbb R}^2, \\
                  u & = 0  \,\,             & \text{in} \,\,    &\partial \Omega.              
\end {split}\right.
$$

Here:

$$ \beta \in (0, +\infty), \,\, 0 \in \partial \Omega, $$
and,

$$ u \in W_0^{1,1}(\Omega), \,\, |x|^{2\beta} e^u \in L^1({\Omega}) \,\, {\rm and} \,\,  0 \leq V \leq b. $$

We can see in [8] a nice formulation of this problem $ (P) $  in the sens of the distributions. This Problem arises from geometrical and physical problems, see for example [1, 3, 21, 25]. The above equation was studied by many authors, with or without  the boundary condition, also for Riemannian surfaces,  see [1-25],  where one can find some existence and compactness results. In [7] we have the following important Theorem ($\beta=0 $):

\smallskip

{\bf Theorem A}{\it (Brezis-Merle [7])}.{\it For $ (u_i)_i $ and $ (V_i)_i $ two sequences of functions relative to $ (P) $ with,
$$ 0 < a \leq V_i \leq b < + \infty $$
then for all compact set $ K $ of $ \Omega $ it holds,
$$ \sup_K u_i \leq c, $$
with $ c $ depending on $ a, b, K $ and $ \Omega $.}

One can find in [7] an interior estimate if we assume $ a=0 $, but we need an assumption on the integral of $ e^{u_i} $, namely, we have ($ \beta=0 $):

\smallskip

{\bf Theorem B}{\it (Brezis-Merle [7])}.{\it For $ (u_i)_i $ and $ (V_i)_i $ two sequences of functions relative to the problem $ (P) $ with,
$$ 0 \leq V_i \leq b < + \infty \,\, {\rm and} \,\, \int_{\Omega} e^{u_i} dy  \leq C, $$
then for all compact set $ K $ of $ \Omega $ it holds;
$$ \sup_K u_i \leq c, $$
with $ c $ depending on $ b, C, K $ and $ \Omega $.}

\smallskip

We look to the uniform boundedness in all $ \bar \Omega $ of the solutions of the Problem $ (P) $. When $ a=0 $, the boundedness of $ \int_{\Omega} e^{u_i} $ is a necessary condition in the problem $ (P) $ as showed in $ [7] $ by the following counterexample ($\beta=0$):

\bigskip

{\bf Theorem C}{\it (Brezis-Merle [7])}.{\it There are two sequences $ (u_i)_i $ and $ (V_i)_i $ of the problem $ (P) $ with, 
$$ 0 \leq V_i \leq b < + \infty \,\, {\rm and} \,\, \int_{\Omega} e^{u_i} dy  \leq C, $$
such that,
$$ \sup_{\Omega}  u_i \to + \infty. $$}

To obtain the two first previous results (Theorems A and B) Brezis and Merle used  an inequality (Theorem 1 of [7]) obtained by an approximation argument and used Fatou's lemma and applied the maximum principle in $ W_0^{1,1}(\Omega) $ which arises from Kato's inequality. Also this weak form of the maximum principle is used to prove the local uniform boundedness result by comparing  a certain function and the Newtonian potential. We refer to [6] for a topic about the weak form of the maximum principle.

When $ \beta = 0 $, the above equation has many properties in the constant and the Lipschitzian cases:

Note that for the problem $ (P) $ ($ \beta =0 $), by using the Pohozaev identity, we can prove that $ \int_{\Omega} e^{u_i} $ is uniformly bounded when $ 0 < a \leq V_i \leq b < +\infty $ and $  ||\nabla V_i||_{L^{\infty}} \leq A $ and $ \Omega $ starshaped, when $ a=0 $ and $ \nabla \log V_i $ is uniformly bounded, we can bound uniformly $ \int_{\Omega} V_i e^{u_i} $  . In [20], Ma-Wei have proved that those results stay true for all open sets not necessarily starshaped in the case $ a>0 $.

\smallskip

In [10] ($ \beta = 0 $) Chen-Li have proved that if $ a=0 $, $ \nabla \log V_i $ is uniformly bounded and $ u_i $ is locally uniformly bounded in $ L^1 $, then  the functions are uniformly bounded near the boundary.

\smallskip

In [10] ($ \beta = 0 $) Chen-Li have proved that if $ a=0 $ and $ \int_{\Omega} e^{u_i} $ is uniformly bounded and $ \nabla \log V_i $ is uniformly bounded, then we have the compactness result directly. Ma-Wei in [20], extend this result in the case where $ a >0 $.

\bigskip

If we assume $ V $ more regular, we can have another type of estimates, a $ \sup + \inf $ type inequalities. It was proved by Shafrir see [23], that, if $ (u_i)_i, (V_i)_i $ are two sequences of functions solutions of the previous equation without assumption on the boundary and, $ 0 < a \leq V_i \leq b < + \infty $, then we have the following interior estimate:
$$ C\left (\dfrac{a}{b} \right ) \sup_K u_i + \inf_{\Omega} u_i \leq c=c(a, b, K, \Omega). $$
One can see in [11] an explicit value of $ C\left (\dfrac{a}{b}\right ) =\sqrt {\dfrac{a}{b}} $. In his proof, Shafrir has used a blow-up function, the Stokes formula and an isoperimetric inequality, see [3]. For Chen-Lin, they have used the blow-up analysis combined with some geometric type inequality for the integral curvature.

\bigskip

Now, if we suppose $ (V_i)_i $ uniformly Lipschitzian with $ A $ the
Lipschitz constant, then, $ C(a/b)=1 $ and $ c=c(a, b, A, K, \Omega)
$, see Brezis-Li-Shafrir [5]. This result was extended for
H\"olderian sequences $ (V_i)_i $ by Chen-Lin, see  [11]. Also, one
can see in [18], an extension of the Brezis-Li-Shafrir result to compact
Riemannian surfaces without boundary. One can see in [19] explicit form,
($ 8 \pi m, m\in {\mathbb N}^* $ exactly), for the numbers in front of
the Dirac masses when the solutions blow-up. Here, the notion of isolated blow-up point is used.

In [9] we have some a priori estimates on the 2 and 3-spheres $ {\mathbb S}_2 $, $ {\mathbb S}_3 $.

Here we give the behavior of the blow-up points on the boundary and a proof of a Brezis-Merle's type Problem when $ \beta \geq 0 $.

We have the following similar problem to Brezis and Merle problem: 

{\bf Problem}. Suppose that $ V_i \to  V $ in $ C^0( \bar \Omega ) $, with, $ 0 \leq V_i $. Also, we consider a sequence of solutions $ (u_i) $ of $ (P) $ relatively to $ (V_i) $ such that,

$$ \int_{\Omega} |x|^{2\beta} e^{u_i} dx \leq C,  $$

is it possible to have:

$$ ||u_i||_{L^{\infty}}\leq C ? $$

Here, we give a caracterization of the behavior of the blow-up points on the boundary and also a proof of the compactness theorem when $ V_i $ are uniformly Lipschitzian. For the behavior of the blow-up points on the boundary, the following condition is enough,

$$ 0 \leq  V_i \leq b, $$

The condition $ V_i \to  V $ in $ C^0(\bar \Omega) $ is not necessary.

But for the proof of the compactness result (for the Brezis-Merle type problem) we assume that:

$$ ||\nabla V_i||_{L^{\infty}}\leq  A. $$

We have the following caracterization of the behavior of the blow-up points on the boundary.

\begin{Theorem}  Assume that $ \max_{\Omega} u_i \to +\infty $, where $ (u_i) $ are solutions of the probleme $ (P) $ with:
 
 $$ \beta \in (0,+\infty), \,\, 0 \leq V_i \leq b,\,\,\, {\rm and } \,\,\, \int_{\Omega}  |x|^{2\beta} e^{u_i} dx \leq C, \,\,\, \forall \,\, i, $$
 
 then;  after passing to a subsequence, there is a finction $ u $,  there is a number $ N \in {\mathbb N} $ and $ N  $ points $ x_1, x_2, \ldots, x_N \in  \partial \Omega $, such that, 
$$ \partial_{\nu} u_i  \to \partial_{\nu} u +\sum_{j=1}^N \alpha_j \delta_{x_j}, \,\,\, \alpha_j \geq 4\pi \,\, {\rm weakly\,\, in \, the \, sens \, of \, measures }. $$
and,
$$ u_i \to u \,\,\, {\rm in }\,\,\, C^1_{loc}(\bar \Omega-\{x_1,\ldots, x_N \}). $$

\end{Theorem} 

 In the following theorem, we have a proof of a global a priori estimate which concern the problem $ (P) $.

\begin{Theorem}Assume that $ (u_i) $ are solutions of $ (P) $ relative to $ (V_i) $ with the following conditions:
$$ \beta \in (0,+\infty),\, 0 \in \partial \Omega, $$
and,
$$ 0 \leq  V_i \leq b, \,\,  ||\nabla V_i||_{L^{\infty}} \leq A,\,\, {\rm and } \,\,\, \int_{\Omega} |x|^{2\beta} e^{u_i} \leq C, $$
We have,
$$  || u_i||_{L^{\infty}} \leq c(b, \beta, A, C, \Omega), $$

\end{Theorem} 

\section{Proof of the theorems} 

\underbar {\it Proof of theorem 1.1:} 

We have, 

$$ u_i\in W^{1,1}_0(\Omega), $$

Since $ \int_{\Omega} |x|^{2\beta} e^{u_i} \leq C $, we have, by the Brezis-Merle result see [7], $ e^{k u_i } \in L^1(\Omega) $ for $ k >2 $ and the elliptic estimates imply that

$$ u_i \in W^{2, k}(\Omega) \cap C^{1, \epsilon}(\bar \Omega). $$

We denote by $ \partial_{\nu} u_i $ the inner normal derivative of $ u_i $. By the maximum principle, $ \partial_{\nu} u_i \geq 0 $.

By the Stokes formula, we obtain 

$$ \int_{\partial \Omega} \partial_{\nu} u_i d\sigma \leq C. $$

Thus, (using the weak convergence in the space of Radon measures), we have the existence of a nonegative Radon measure $ \mu $ such that

$$ \int_{\partial \Omega} \partial_{\nu} u_i d\sigma \leq C, $$

Without loss of generality, we can assume that $ \partial_{\nu} u_i >0 $. Thus, (using the weak convergence in the space of Radon measures), we have the existence of a positive Radon measure $ \mu $ such that,

$$ \int_{\partial \Omega} \partial_{\nu} u_i \phi  d\sigma \to \mu(\phi), \,\,\, \forall \,\,\, \phi \in C^0(\partial \Omega). $$

We take an $ x_0 \in \partial \Omega $ such that, $ \mu({x_0}) < 4\pi $. Without loss of generality, we can assume that the following curve, $ B(x_0, \epsilon) \cap \partial \Omega := I_{\epsilon} $ is an interval.

We take a test function:

$$ \begin{cases}
    
\eta_{\epsilon} \equiv 1,\,\,\,  {\rm on } \,\,\, I_{\epsilon}, \,\,\, 0 < \epsilon < \delta/2,\\

\eta_{\epsilon} \equiv 0,\,\,\, {\rm outside} \,\,\,I_{2\epsilon }, \\

0 \leq \eta_{\epsilon} \leq 1, \\

||\nabla \eta_{\epsilon}||_{L^{\infty}(I_{2\epsilon})} \leq \dfrac{C_0(\Omega, x_0)}{\epsilon}.

\end{cases} $$

We take a $\tilde \eta_{\epsilon} $ such that,

$$  \left \{ \begin {split} 
       -\Delta \tilde \eta_{\epsilon} &= 0 && \text{in} \,\,\Omega \\
            \tilde\eta_{\epsilon} &=  \eta_{\epsilon} && \text{on} \,\, \partial \Omega.              
\end {split}\right.
$$

{\bf Remark:} We use the following steps in the construction of $ \tilde \eta_{\epsilon} $:

We take a cutoff function $ \eta_{0} $ in $ B(0, 2) $ or $ B(x_0, 2) $:

1- We set $ \eta_{\epsilon}(x)= \eta_0(|x-x_0|/\epsilon) $ in the case of the unit disk it is sufficient.

2- Or, in the general case: we use a chart $ (f, \tilde \Omega) $ with $ f(0)=x_0 $ and we take $ \mu_{\epsilon}(x)= \eta_0 ( f( |x|/ \epsilon)) $ to have  connected  sets $ I_{\epsilon} $ and we take $ \eta_{\epsilon}(y)= \mu_{\epsilon}(f^{-1}(y))$. Because $ f, f^{-1} $ are Lipschitz, $ |f(x)-x_0| \leq k_ 2|x|\leq 1 $ for $ |x| \leq 1/k_2 $ and $ |f(x)-x_0| \geq k_ 1|x|\geq 2 $ for $ |x| \geq 2/k_1>1/k_2 $, the support  of $ \eta $ is in $ I_{(2/k_1)\epsilon} $.

$$ \begin{cases}
    
\eta_{\epsilon} \equiv 1,\,\,\,  {\rm on } \,\,\,  f(I_{(1/k_2)\epsilon}), \,\,\, 0 < \epsilon < \delta/2,\\

\eta_{\epsilon} \equiv 0,\,\,\, {\rm outside} \,\,\, f(I_{(2/k_1)\epsilon }), \\

0 \leq \eta_{\epsilon} \leq 1, \\

||\nabla \eta_{\epsilon}||_{L^{\infty}(I_{(2/k_1)\epsilon})} \leq \dfrac{C_0(\Omega, x_0)}{\epsilon}.

\end{cases} $$

3- Also, we can take: $ \mu_{\epsilon}(x)= \eta_0(x/\epsilon) $ and $ \eta_{\epsilon}(y)= \mu_{\epsilon}(f^{-1}(y)) $, we extend it by $ 0 $ outside $ f(B_1(0)) $.  We have $ f(B_1(0)) = D_1(x_0) $, $ f (B_{\epsilon}(0))= D_{\epsilon}(x_0) $ and $ f(B_{\epsilon}^+)= D_{\epsilon}^+(x_0) $ with $ f $ and $ f^{-1} $ smooth diffeomorphism.

$$ \begin{cases}
    
\eta_{\epsilon} \equiv 1,\,\,\,  {\rm on \, a \, the \, connected \, set } \,\,\,  J_{\epsilon} =f(I_{\epsilon}), \,\,\, 0 < \epsilon < \delta/2,\\

\eta_{\epsilon} \equiv 0,\,\,\, {\rm outside} \,\,\, J'_{\epsilon}=f(I_{2\epsilon }), \\

0 \leq \eta_{\epsilon} \leq 1, \\

||\nabla \eta_{\epsilon}||_{L^{\infty}(J'_{\epsilon})} \leq \dfrac{C_0(\Omega, x_0)}{\epsilon}.

\end{cases} $$

And, $ H_1(J'_{\epsilon}) \leq C_1 H_1(I_{2\epsilon}) = C_1 4\epsilon $, since $ f $ is Lipschitz. Here $ H_1 $ is the Hausdorff measure.

 We solve the Dirichlet Problem:

\begin{displaymath}  \left \{ \begin {split} 
      \Delta \bar \eta_{\epsilon}  & = \Delta \eta_{\epsilon}              \,\, &&\text{in} \!\!&&\Omega \subset {\mathbb R}^2, \\
                  \bar \eta_{\epsilon} & = 0   \,\,             && \text{in} \!\!&&\partial \Omega.               
\end {split}\right.
\end{displaymath}

and finaly we set $ \tilde \eta_{\epsilon} =-\bar \eta_{\epsilon} + \eta_{\epsilon} $. Also, by the maximum principle and the elliptic estimates we have :

$$ ||\nabla \tilde \eta_{\epsilon}||_{L^{\infty}} \leq C(|| \eta_{\epsilon}||_{L^{\infty}} +||\nabla \eta_{\epsilon}||_{L^{\infty}} + ||\Delta \eta_{\epsilon}||_{L^{\infty}}) \leq \dfrac{C_1}{\epsilon^2}, $$

with $ C_1 $ depends on $ \Omega $.

We use the following estimate, see [4,8,14,25]

$$ ||\nabla u_i||_{L^q} \leq C_q, \,\,\forall \,\, i\,\, {\rm and  }  \,\, 1< q < 2. $$

We deduce from the last estimate that, $ (u_i) $ converge weakly in $ W_0^{1, q}(\Omega) $, almost everywhere to a function $ u \geq 0 $ and $ \int_{\Omega} |x|^{2\beta} e^u < + \infty $ (by Fatou lemma). Also, $ V_i $ weakly converge to a nonnegative function $ V $ in $ L^{\infty} $. The function $ u $ is in $ W_0^{1, q}(\Omega) $ solution of :

$$  \left \{ \begin {split} 
       -\Delta u &= |x|^{2\beta}Ve^{u} \in L^1(\Omega) && {\rm in } \,\, \Omega \\
                   u&= 0 && \text{on} \,\,\partial \Omega,              
\end {split}\right.
$$

As in  the corollary 1 of Brezis-Merle result, see [7],   we have $ e^{k u }\in L^1(\Omega), k >1 $. By the elliptic estimates, we have $ u \in C^1(\bar \Omega) $.

We can write,

 \be -\Delta ((u_i-u) \tilde \eta_{\epsilon})= |x|^{2\beta}(V_i e^{u_i} -Ve^u)\tilde \eta_{\epsilon} +2<\nabla (u_i- u)| \nabla \tilde \eta_{\epsilon}> .\label{(1)}\ee

We use the interior esimate of Brezis-Merle, see [7],

\underbar {\it Step 1:} Estimate of the integral of the first term of the right hand side of $ (1) $.

We use the Green formula between $ \tilde \eta_{\epsilon} $ and $ u $, we obtain,

 \be \int_{\Omega} |x|^{2\beta}Ve^u \tilde \eta_{\epsilon} dx =\int_{\partial \Omega} \partial_{\nu} u \eta_{\epsilon} \leq 4\epsilon ||\partial_{\nu}u||_{L^{\infty}}= C \epsilon \label{(2)}\ee

We have,

$$  \left \{ \begin {split} 
       -\Delta u_i &= |x|^{2\beta}V_ie^{u_i} && {\rm in } \,\, \Omega \\
                   u_i&= 0 && \text{on} \,\,\partial \Omega,              
\end {split}\right.
$$

We use the Green formula between $ u_i $ and $ \tilde \eta_{\epsilon} $ to have:

\be \int_{\Omega} |x|^{2\beta}V_i e^{u_i} \tilde \eta_{\epsilon} dx = \int_{\partial \Omega} \partial_{\nu} u_i \eta_{\epsilon} d\sigma \to \mu(\eta_{\epsilon}) \leq \mu(I_{2\epsilon}) \leq 4  \pi - \epsilon_0, \,\,\, \epsilon_0 >0 \label{(3)}\ee

From $ (2) $ and $ (3) $ we have for all $ \epsilon >0 $ there is $ i_0 =i_0(\epsilon) $ such that, for $ i \geq i_0 $,

\be \int_{\Omega} |x|^{2\beta}|(V_ie^{u_i}-Ve^u) \tilde \eta_{\epsilon}| dx \leq 4 \pi -\epsilon_0 + C \epsilon \label{(4)} \ee

\underbar { Step 2:} Estimate of integral of the second term of the right hand side of $ (\ref{(1)}) $.

Let $ \Sigma_{\epsilon} = \{ x \in \Omega, d(x, \partial \Omega) = \epsilon^3 \} $ and $ \Omega_{\epsilon^3} = \{ x \in \Omega, d(x, \partial \Omega) \geq \epsilon^3 \} $, $ \epsilon > 0 $. Then, for $ \epsilon $ small enough, $ \Sigma_{\epsilon} $ is hypersurface.

The measure of $ \Omega-\Omega_{\epsilon^3} $ is $ k_2\epsilon^3 \leq \mu_L (\Omega-\Omega_{\epsilon^3}) \leq k_1 \epsilon^3 $.

{\bf Remark}: for the unit ball $ {\bar B(0,1)} $, our new manifold is $ {\bar B(0, 1-\epsilon^3)} $.

(Proof of this fact; let's consider $ d(x, \partial \Omega) = d(x, z_0), z_0 \in \partial \Omega $, this imply that $ (d(x,z_0))^2 \leq (d(x, z))^2 $ for all $ z \in \partial \Omega $ which it is equivalent to $ (z-z_0) \cdot (2x-z-z_0) \leq 0 $ for all $ z \in \partial \Omega $, let's consider a chart around $ z_0 $ and $ \gamma (t) $ a curve in $ \partial \Omega $, we have;

$ (\gamma (t)-\gamma(t_0) \cdot (2x-\gamma(t)-\gamma(t_0)) \leq 0 $ if we divide by $ (t-t_0) $ (with the sign and tend $ t $ to $ t_0 $), we have $ \gamma'(t_0) \cdot (x-\gamma(t_0)) = 0 $, this imply that $ x= z_0-s \nu_0 $ where $ \nu_0 $ is the outward normal of $ \partial \Omega $ at $ z_0 $))

With this fact, we can say that $ S= \{ x, d(x, \partial \Omega) \leq \epsilon \}= \{ x= z_0- s \nu_{z_0}, z_0 \in \partial \Omega, \,\, -\epsilon \leq s \leq \epsilon \} $. It  is sufficient to work on  $ \partial \Omega $. Let's consider a charts $ (z, D=B(z, 4 \epsilon_z), \gamma_z) $ with $ z \in \partial \Omega $ such that $ \cup_z B(z, \epsilon_z) $ is  cover of $ \partial \Omega $ .  One can extract a finite cover $ (B(z_k, \epsilon_k)), k =1, ..., m $, by the area formula the measure of $ S \cap B(z_k, \epsilon_k) $ is less than a $ k\epsilon $ (a $ \epsilon $-rectangle).  For the reverse inequality, it is sufficient to consider one chart around one point of the boundary).

We write,

$$ \int_{\Omega} |<\nabla ( u_i -u)|\nabla \tilde \eta_{\epsilon} > | dx =
\int_{\Omega_{\epsilon^3}} |<\nabla (u_i -u)|\nabla \tilde \eta_{\epsilon} > | dx + $$

\be + \int_{\Omega - \Omega_{\epsilon^3}} <\nabla (u_i-u)|\nabla \tilde \eta_{\epsilon} >| dx. \label{(5)} \ee

\underbar {\it Step 2.1:} Estimate of $ \int_{\Omega - \Omega_{\epsilon^3}} |<\nabla (u_i-u)|\nabla \tilde \eta_{\epsilon} >| dx $.

First, we know from the elliptic estimates that  $ ||\nabla \tilde \eta_{\epsilon}||_{L^{\infty}} \leq C_1 /\epsilon^2 $, $ C_1 $ depends on $ \Omega $.

We know that $ (|\nabla u_i|)_i $ is bounded in $ L^q, 1< q < 2 $, we can extract  from this sequence a subsequence which converge weakly to $ h \in L^q $. But, we know that we have locally the uniform convergence to $ |\nabla u| $ (by Brezis-Merle theorem), then, $ h= |\nabla u| $ a.e. Let $ q' $ be the conjugate of $ q $.

We have, $  \forall f \in L^{q'}(\Omega)$

$$ \int_{\Omega} |\nabla u_i| f dx \to \int_{\Omega} |\nabla u| f dx $$

If we take $ f= 1_{\Omega-\Omega_{\epsilon^3}} $, we have:

$$ {\rm for } \,\, \epsilon >0 \,\, \exists \,\, i_1 = i_1(\epsilon) \in {\mathbb N}, \,\,\, i \geq  i_1,  \,\, \int_{\Omega-\Omega_{\epsilon^3}} |\nabla u_i| \leq \int_{\Omega-\Omega_{\epsilon^3}} |\nabla u| + \epsilon^3. $$

Then, for $ i \geq i_1(\epsilon) $,

$$ \int_{\Omega-\Omega_{\epsilon^3}} |\nabla u_i| \leq mes(\Omega-\Omega_{\epsilon^3}) ||\nabla u||_{L^{\infty}} + \epsilon^3 = \epsilon^3(k_1 ||\nabla u||_{L^{\infty}} + 1). $$

Thus, we obtain,

\be \int_{\Omega - \Omega_{\epsilon^3}} |<\nabla (u_i-u)|\nabla \tilde \eta_{\epsilon} >| dx \leq  \epsilon C_1(2 k_1 ||\nabla u||_{L^{\infty}} + 1) \label{(6)}\ee

The constant $ C_1 $ does  not depend on $ \epsilon $ but on $ \Omega $.

\underbar {\it Step 2.2:} Estimate of $ \int_{\Omega_{\epsilon^3}} |<\nabla (u_i-u)|\nabla \tilde \eta_{\epsilon} >| dx $.

We know that, $ \Omega_{\epsilon} \subset \subset \Omega $, and ( because of Brezis-Merle's interior estimates) $ u_i \to u $ in $ C^1(\Omega_{\epsilon^3}) $. We have,

$$ ||\nabla (u_i-u)||_{L^{\infty}(\Omega_{\epsilon^3})} \leq \epsilon^3,\, {\rm for } \,\, i \geq i_3 = i_3(\epsilon). $$

We write,
 
$$ \int_{\Omega_{\epsilon^3}} |<\nabla (u_i-u)|\nabla \tilde \eta_{\epsilon} >| dx \leq ||\nabla (u_i-u)||_{L^{\infty}(\Omega_{\epsilon^3})} ||\nabla \tilde \eta_{ \epsilon}||_{L^{\infty}} \leq C_1 \epsilon \,\, {\rm for } \,\, i \geq i_3, $$

For $ \epsilon >0 $, we have for $ i \in {\mathbb N} $, $ i \geq \max \{i_1, i_2, i_3 \} $,

\be \int_{\Omega} |<\nabla (u_i-u)|\nabla \tilde \eta_{\epsilon} >| dx \leq \epsilon C_1(2 k_1 ||\nabla u||_{L^{\infty}} + 2) \label{(7)}\ee

From $ (4) $ and $ (7) $, we have, for $ \epsilon >0 $, there is $ i_3= i_3(\epsilon) \in {\mathbb N}, i_3 = \max \{ i_0, i_1, i_2 \} $ such that,

\be \int_{\Omega} |\Delta [(u_i-u)\tilde \eta_{\epsilon}]|dx \leq 4 \pi-\epsilon_0+  \epsilon 2 C_1(2 k_1 ||\nabla u||_{L^{\infty}} + 2 + C) \label{(8)}\ee

We choose $ \epsilon >0 $ small enough to have a good estimate of  $ (1) $.

Indeed, we have:

$$ \left \{ \begin {split}
 -\Delta [(u_i-u) \tilde \eta_{\epsilon}] &= g_{i,\epsilon} && {\rm in } \,\,\, \Omega,\\
 (u_i-u) \tilde \eta_{\epsilon} &= 0 && {\rm on } \,\,\, \partial \Omega.
 \end {split}  \right . $$

with $ ||g_{i, \epsilon} ||_{L^1(\Omega)} \leq 4 \pi -\epsilon_0. $

We can use Theorem 1 of [7] to conclude that there is $ q >1 $ such that:

$$ \int_{V_{\epsilon}(x_0)} e^{q|u_i-u|} dx \leq \int_{\Omega} e^{q|u_i -u| \tilde \eta_{\epsilon}} dx \leq C(\epsilon,\Omega). $$
 
where, $ V_{\epsilon}(x_0) $ is a neighberhood of $ x_0 $ in $ \bar \Omega $.

Thus, for each $ x_0 \in \partial \Omega - \{ \bar x_1,\ldots, \bar x_m \} $ there is $ \epsilon_{x_0} >0, q_1 > 1 $ such that:

$$ \int_{B(x_0, \epsilon_{x_0})} e^{q_1 u_i} dx \leq C, \,\,\, \forall \,\,\, i. $$

Now, we consider a cutoff function $ \eta \in C^{\infty}({\mathbb R}^2) $ such that

$$ \eta \equiv 1 \,\,\, {\rm on } \,\,\, B(x_0, \epsilon_{x_0}/2) \,\,\, {\rm and } \,\,\, \eta \equiv 0 \,\,\, {\rm on } \,\,\, {\mathbb R}^2 -B(x_0, 2\epsilon_{x_0}/3). $$

We write

$$ -\Delta (u_i \eta) = |x|^{2\beta}V_i e^{u_i} \eta - 2 \nabla u_i \cdot \nabla \eta  - u_i \Delta \eta. $$

By the elliptic estimates (see [15]) $ (u_i)_i $ is uniformly bounded in $ W^{2, q_1}(V_{\epsilon}(x_0)) $ and also, in $ C^1(V_{\epsilon}(x_0)) $. Finaly, we have, for some $ \epsilon > 0 $ small enough,

$$ || u_i||_{C^{1,\theta}[B(x_0, \epsilon)]} \leq c_3 \,\,\, \forall \,\,\, i. $$

\underbar {\it Proof of theorem 1.2:} 

Without loss of generality, we can assume that $ 0 $ is a blow-up point. Since the boundary is an analytic curve $ \gamma(t) $, there is a neighborhood of  $ 0 $ such that the curve $ \gamma $ can be extend to a holomorphic map such that $ \gamma'(0) \not = 0 $ (series) and by the inverse mapping one can assume that this map is univalent around $ 0 $. In the case when the boundary is a simple Jordan curve the domain is simply connected, see [24]. In the case that the domains has a finite number of holes it is conformally equivalent to a disk with a finite number of disks removed, see [17]. Here we consider a general domain. Without loss of generality one can assume that $ \gamma (B_1^+) \subset \Omega $ and also $ \gamma (B_1^-) \subset (\bar \Omega)^c $ and $ \gamma (-1,1) \subset \partial \Omega $ and $ \gamma $ is univalent. This means that $ (B_1, \gamma) $ is a local chart around $ 0 $ for $ \Omega $ and $ \gamma $ univalent. (This fact holds if we assume that we have an analytic domain, in the sense of Hofmann see [16], (below a graph of an analytic function), we have necessary the condition $ \partial \bar \Omega = \partial \Omega $ and the graph is analytic, in this case $ \gamma (t)= (t, \phi(t)) $ with $ \phi $ real analytic and an example of this fact is the unit disk  around the point $ (0,1) $ for example).

By this conformal transformation, we can assume that $ \Omega =B_1^+ $, the half ball, and $ \partial^+ B_1^+ $ is the exterior part, a part which not contain $ 0 $ and on which  $ u_i $ converge in the $ C^1 $ norm to $ u $. Let us consider $ B_{\epsilon}^+ $, the half ball with radius $ \epsilon >0 $. Also, one can consider a $ C^1 $ domain (a rectangle between two half disks) and by charts its image is a $ C^1 $ domain).

We know that:

$$ u_i \in C^{2, \epsilon}(\bar \Omega). $$ 

Thus we can use integrations by parts (Gauss-Green-Riemann-Stokes formula). The second Pohozaev identity applied around the blow-up $ 0 $ see for example [2, 20, 22] gives :

$$ 2(1+ \beta ) \int_{B_{\epsilon}^+} |x|^{2 \beta} V_ie^{u_i} dx + \int_{B_{\epsilon}^+}  < x |\nabla V_i > |x|^{2 \beta} V_ie^{u_i} dx + \int_{\partial B_{\epsilon}^+} < \nu |x >|x|^{2 \beta}V_ie^{u_i} d\sigma= $$

\be = \int_{\partial^+ B_{\epsilon}^+}  g(\nabla u_i)d\sigma, \label{(10)}\ee

with,

$$ g(\nabla u_i)=< \nu |\nabla u_i >< x |\nabla u_i >-< x |\nu > \dfrac{|\nabla u_i|^2}{2}. $$

also,

$$ 2(1+ \beta ) \int_{B_{\epsilon}^+} |x|^{2 \beta} Ve^{u} dx + \int_{B_{\epsilon}^+}  < x |\nabla V> |x|^{2 \beta} Ve^{u} dx + \int_{\partial B_{\epsilon}^+} < \nu |x >|x|^{2 \beta}Ve^{u} d\sigma = $$
\be = \int_{\partial^+ B_{\epsilon}^+}  g(\nabla u)d\sigma, \label{(10)}\ee

Thus, 

$$ 2(1+ \beta ) \int_{B_{\epsilon}^+} |x|^{2 \beta} V_ie^{u_i} dx-2(1+ \beta ) \int_{B_{\epsilon}^+} |x|^{2 \beta} Ve^{u} dx + $$

$$ + \int_{B_{\epsilon}^+}  < x |\nabla V_i > |x|^{2 \beta} V_ie^{u_i} dx - \int_{B_{\epsilon}^+}  < x |\nabla V> |x|^{2 \beta} Ve^{u} dx+ o(1) = $$

$$  = \int_{\partial^+ B_{\epsilon}^+}  g(\nabla u_i)-g(\nabla u) d\sigma = o(1), $$

First, we tend $ i  $ to infinity after $ \epsilon  $ to 0, we obtain:

\be \lim_{ \epsilon \to 0} \lim_{ i\to + \infty} 2(1+ \beta ) \int_{B_{\epsilon}^+} |x|^{2 \beta} V_ie^{u_i} dx = 0, \label{(12)}\ee

however

\be \int_{ \gamma ( B_{\epsilon}^+)} |x|^{2 \beta}V_i e^{u_i} dx =\int_{\partial  \gamma ( B_{\epsilon}^+)} \partial_{\nu} u_i d\sigma = \alpha_1+ O(\epsilon)+o(1) >0. \ee

which is a contradiction.

Here we used a theorem of Hofmann see [16], which gives the fact that $ \gamma (B_{\epsilon}^+) $ is a Lipschitz domain. Also, we can see that $ \gamma ((-\epsilon, \epsilon)) $ and $ \gamma (\partial^+ B_{\epsilon}^+) $ are submanifolds.  

We start with a Lipschitz domain $ B_{\epsilon}^+ $ because it is convex and by the univalent and conformal map $ \gamma $ the image of this domain $ \gamma (B_{\epsilon}^+) $ is a Lipschitz domain and thus we can apply the integration by part and here we know the explicit formula of the unit outward normal it is the usual unit outward normal (normal to the tangent space of the boundary which we know explicitly because we have two submanifolds).

In the case of the disk $ D = \Omega $, it is sufficient to consider $ B(0,\epsilon)  \cap D $ which is a Lipschitz domain because it is convex (and not necessarily $ \gamma (B_{\epsilon}^+) $).

There is a version of the integration by part which is the Green-Riemann formula in dimension 2 on a domain $ \Omega $. This formula holds if we assume that there is a finite number of points $ y_1,..., y_m $ such that $ \partial \Omega - (y_1,..., y_m) $ is a $ C^1 $ manifold and for $ C^1 $ tests functions, see [2], for the Gauss-Green-Riemann-Stokes formula, for $ C^1 $ domains with singular points (here a finite number of singular points).

{\bf Remark 1:} Note that a monograph of Droniou contain a proof of all fact about Sobolev spaces (with Strong Lipschitz property) with only weak Lipschitz property (Lipschitz-Charts), we start with Strong Lipschitz property and by $ \gamma $ we have weak Lipschtz property.

{\bf Remarks about the conformal map :}

1-It sufficient to prove  that $ \gamma_1((-\epsilon,\epsilon))=
\partial \Omega \cap \tilde \gamma_1(B_{\epsilon})= \partial \Omega \cap \tilde \gamma_1(B_{\epsilon})\cap \{ |abscissa|< \epsilon \} $, for $ \epsilon >0 $ small enough. Where $ \tilde \gamma_1 $ is the holomorphic extension of $ \gamma_1(t)=t+i\phi(t) $. For this, we argue by contradiction, we have for $ z_{\epsilon} \in B_{\epsilon} $, $\tilde \gamma_1(z_{\epsilon}) =(t_{\epsilon},\phi (t_{\epsilon})) $ for $ |t_{\epsilon}|\geq \epsilon $. Because $ \tilde \gamma_1 $ is injective on $ B_1 $ and $ \tilde \gamma_1 =\gamma_1=t+i\phi(t) $ on the real axis, we have necessirely $ |t_{\epsilon}|\geq 1 $. But, by continuity $ |\tilde \gamma_1(z_{\epsilon})|\to 0 $ because $ z_{\epsilon} \to 0 $. And, we use the fact that $ |\tilde \gamma_1(z_{\epsilon})|=|(t_{\epsilon},\phi(t_{\epsilon}))|\geq |t_{\epsilon}|\geq 1 $, to have a contradiction.) (This means  that for a small radius when the graph go out from the ball, it never retruns to the ball). (This fact imply that, when we have a curve which cut $\partial \Omega $ in $ \tilde \gamma_1(B_{\epsilon}) $ then the point have an abscissa such that $ |abscissa|< \epsilon $. This fact (by a contradiction with the fact $ \partial \bar \Omega = \partial \Omega $ and consider paths), imply that the image of the upper part of the ball is one side of the curve and the image of the lower part is in the other side of the curve.

2- Also, we can consider directly the coordinate $ T $ and change the function $ u(z)\to u(T) $ by $ z=T/a_1+x_0 $.

Set: $ \psi(\lambda_1,\lambda_2)\to M \in \Omega $ such that, $ \overrightarrow{x_0M}=\lambda_1 i_1'+\lambda_2j_1' $ where $ (i_1',j_1') $ is a basis such that $ i_1'=e^{-i\theta} i_1, j_1'=e^{-i\theta} j_1 $. And, $ \phi(x_1,x_2)\to M $ such that, $ \overrightarrow{OM}=x_1i_1+x_2i_2 $ the canonical basis $ (i_1,j_1) $. Then, we have two charts $ \phi $ and $ \psi $ and the complex affix $ T_M=\lambda_1+i\lambda_2 $ and $ z_M=x_1+ix_2 $ are such that  (transition map):

$$ T_M/e^{i\theta}+x_0=z_M=\phi^{-1} o \psi(\lambda_1,\lambda_2), $$

We have:

$$ \partial_{\lambda_1}=\cos \theta \partial_{x_1}+\sin \theta \partial_{x_2}, $$

$$ \partial_{\lambda_2}=-\sin \theta \partial_{x_1}+\cos \theta \partial_{x_2}, $$

Thus, the metric in the chart $\psi $ or coordinates $ (\lambda_1,\lambda_2) $ is : $ g_{ij}^{\lambda}=\delta_{ij} $ and the Laplacian in the two charts, $ \psi $ and $\phi $ are the usual Laplacian $ \partial_{\lambda_1\lambda_1}+\partial_{\lambda_2\lambda_2} $.

We write:

$$ \Delta u(M)=\Delta_{\lambda} (u o\psi(\lambda_1,\lambda_2)) $$

And then we apply the conformal map $\tilde \gamma_1 $ which send the affix $ T_M $, $ M $ in a neighborhood of $ x_0\in \partial \Omega $ to $ B_{\epsilon} $ with the fact that send $ T_M, M\in \partial \Omega $ to the real axis $ (-\epsilon,\epsilon) $ and the other parts of $ \Omega $ and $ \bar \Omega^c $.

\smallskip

We have: $ \psi:(\lambda_1,\lambda_2) \to M $ and $ \tilde \gamma_1:(\lambda_1,\lambda_2)\to (\mu_1,\mu_2) $. Then the map considerded is: $ \psi o \tilde \gamma_1^{-1} $. It is a chart. Indeed the first chart is: $ \psi o g^{-1} $ where $ g $ is $ g:(\lambda_1,\lambda_2)\to [\lambda_1,\lambda_2+\phi_0(\lambda_1)] $ with $ \phi_0 $ the "graph map" as in the defiinition of $ \partial \Omega $.

Then the transition map is:

$$ g o\tilde \gamma_1^{-1} $$

is smooth from a neighborhood of $ 0 $ to a neighborhood of  $ 0 $ and send the upper part of a ball (lower part) to an upper part of a domain (lower part). And then $ \psi o\tilde \gamma_1^{-1} $ is a chart with the property that $ \tilde \gamma_1 $ is conformal.

Remark that: $ \psi $ and $ \phi $ are charts for  $ \Omega $ and $ \psi $ is almost chart for $ \partial \Omega $. But $ \psi o \tilde \gamma_1^{-1} $ and $ \psi o g^{-1} $ are the "real" charts of the boundary.

\smallskip

3-We can remark that a definition of $ C^k, k \geq 1 $ domain, is equivalent to a definition of a submanifold with the condition $ \partial \bar \Omega = \partial \Omega $ or $ {\dot{\bar \Omega }} = \Omega $.

{\bf Remark 2: about a variational formulation.}

we consider a solutions in the sense of distribution. By the same argument (in the proof of the maximum principle obtained by Kato's inequality $ W^{1,1}_0 $ is sufficient), see [6], we prove that the solutions are in the sense $ C^2_0 $ of Agmon, see [1]. Also, we have corollary 1 of [7]. We use Agmon's regularity theorem. We return to the usual variational formulation in $ W^{1,2}_0 $ and thus we have the estimate of [8] or by Stampacchia duality theorem in $ W^{1,2}_0 $.

\end{document}